\documentclass[11pt,oneside]{amsart}

\usepackage{amssymb,amsthm,amsmath}
\usepackage{pstricks}
\usepackage{pst-node}
\usepackage{pst-tree}
\usepackage{graphicx}
\usepackage{pstcol}
\usepackage{pst-plot}
\usepackage{subfigure}
\usepackage{color}
\definecolor{darkblue}{rgb}{0, 0, .4}
\usepackage[bookmarks,breaklinks]{hyperref}
\hypersetup{
	colorlinks=true,
	linkcolor=darkblue,
	anchorcolor=darkblue,
	citecolor=darkblue,
	pagecolor=darkblue,
	urlcolor=darkblue,
	pdfpagemode=UseThumbs,
	pdftitle={Enumeration schemes for restricted permutations},
	pdfsubject={Enumerative Combinatorics},
	pdfauthor={Vincent Vatter},
	pdfkeywords={enumeration scheme, forbidden subsequence, generating tree, insertion encoding, restricted permutation, simple permutation, Wilfian formula}
}
	
\newtheorem{theorem}{Theorem}[section]
\newtheorem{proposition}[theorem]{Proposition}

\newtheorem{definition}[theorem]{Definition}

\newtheorem{question}[theorem]{Question}

\newcounter{todocounter}

\newcommand{\st}{\operatorname{st}}
\newcommand{\strongcomp}{\operatorname{sc}}
\newcommand{\Av}{\operatorname{Av}}
\newcommand{\Binf}{\|B\|_\infty}
\newcommand{\J}{\mathcal{J}}
\newcommand{\G}{\mathcal{G}}
\newcommand{\C}{\mathcal{C}}

\newcommand{\esp}{ES$^+$}
\renewcommand{\d}{d}
\renewcommand{\sb}{\mathcal{SB}}

\renewcommand{\a}{{\bf a}}
\renewcommand{\b}{{\bf b}}
\newcommand{\g}{{\bf g}}
\newcommand{\h}{{\bf h}}
\newcommand{\x}{{\bf x}}
\newcommand{\y}{{\bf y}}
\newcommand{\case}[4]{\left\{\begin{array}{ll}#1&\mbox{#2}\\#3&\mbox{#4}\end{array}\right.}
\newcommand{\ie}[1]{\mathbf{#1}}
\newcommand{\OEISlink}[1]{\href{http://www.research.att.com/projects/OEIS?Anum=#1}{#1}}
\newcommand{\OEISref}{\href{http://www.research.att.com/\~njas/sequences/}{OEIS}~\cite{OEIS}}
\newcommand{\OEIS}[1]{(Sequence \OEISlink{#1} in the \OEISref.)}
\newcommand{\basisnodeone}[2]{\begin{minipage}[t]{1in}
\begin{center}
$#1$\\
\begin{tabular}{rcl}
$\g\!\!\!\!\!$&$\not\ge$&$\!\!\!\!\!#2$
\end{tabular}
\end{center}
\end{minipage}}
\newcommand{\basisnodetwo}[3]{\begin{minipage}[t]{1in}
\begin{center}
$#1$\\
\begin{tabular}{rcl}
$\g\!\!\!\!\!$&$\not\ge$&$\!\!\!\!\!#2$,\\
&&$\!\!\!\!\!#3$
\end{tabular}
\end{center}
\end{minipage}}
\newcommand{\basisnodethree}[4]{\begin{minipage}[t]{1in}
\begin{center}
$#1$\\
\begin{tabular}{rcl}
$\g\!\!\!\!\!$&$\not\ge$&$\!\!\!\!\!#2$,\\
&&$\!\!\!\!\!#3$,\\
&&$\!\!\!\!\!#4$
\end{tabular}
\end{center}
\end{minipage}}

\begin{document}

\title{Enumeration schemes for restricted permutations}

\author[Vincent Vatter]{Vincent Vatter}
\address{Department of Mathematics, Rutgers University, New Brunswick, New Jersey}
\email{\href{mailto:vatter@math.rutgers.edu}{\texttt{vatter@math.rutgers.edu}}}
\urladdr{\url{http://math.rutgers.edu/\~vatter/}}
\thanks{Partially supported by an award from DIMACS and an NSF VIGRE grant to the Rutgers University Department of Mathematics.}

\date{\today}
\subjclass[2000]{05A05, 05A15, 68Q20}
\keywords{forbidden subsequence, generating tree, insertion encoding, restricted permutation, simple permutation, Wilfian formula, enumeration scheme}

\begin{abstract}
Zeilberger's enumeration schemes can be used to completely automate the enumeration of many permutation classes.  We extend his enumeration schemes so that they apply to many more permutation classes and describe the Maple package {\sc WilfPlus}, which implements this process.  We also compare enumeration schemes to three other systematic enumeration techniques: generating trees, substitution decompositions, and the insertion encoding.
\end{abstract}

\maketitle

\section{Introduction}\label{wp-intro}

The enumeration of permutation classes, whose ancestry can be traced back to at least 1915 (MacMahon~\cite{m:ca}), has frequently been accomplished by beautiful arguments utilizing such diverse objects as Young tableaux, Dyck paths, and planar maps, to name only a few.  Our concern herein is not with attractive proofs, but rather with systematic methods for solving the enumeration problem.  We adopt a strict definition of systematic, insisting that the computations can be performed without any human interaction whatsoever.  For the definition of enumeration, we follow Wilf~\cite{wilf:formula} and insist only on a polynomial time (in $n$) algorithm to compute the number of length $n$ permutations in the class.  We refer to such an algorithm as a {\it Wilfian formula\/}.  To date, four techniques with wide applicability have been introduced which satisfy these goals:
\begin{itemize}
\item generating trees, 
\item enumeration schemes, 
\item substitution decompositions, 
\item the insertion encoding.
\end{itemize}

The major aim of this paper, carried out in Section~\ref{wp-wp}, is to extend the method of enumeration schemes so that it can enumerate a wider variety of permutation classes and describe the Maple package {\sc WilfPlus}, which can rigorously and automatically find these extended schemes.  Before that, we briefly examine the other methods in Sections~\ref{wp-gt}--\ref{wp-simple} and review enumeration schemes in Section~\ref{wp-wilf}.  Section~\ref{wp-nonex} contains examples of classes which lie beyond the reach of even our more powerful enumeration schemes, while Section~\ref{wp-ex} gives numerous examples which can be handled.  First we describe permutation classes.

Two sequences of natural numbers are said to be {\it order isomorphic\/} if they have the same pairwise comparisons, so $9,1,6,7,2$ is order isomorphic to $5,1,3,4,2$.  Every sequence $w$ of natural numbers without repetition is order isomorphic to a unique permutation that we denote by $\st(w)$, so $\st(9,1,6,7,2)=5,1,3,4,2$, which we shorten to $51342$.  We say that $\st(w)$ is the {\it standardization\/} of $w$.  We further say that the permutation $\pi$ {\it contains\/} the permutation $\beta$ if $\pi$ contains a subsequence that is order isomorphic to $\beta$, and in this case we write $\beta\le\pi$. For example, $391867452$ contains $51342$, as can be seen by considering the subsequence $91672$.  A permutation is said to {\it avoid\/} another if it does not contain it.

A {\it permutation class\/} is a lower order ideal in the containment ordering, meaning that if $\pi$ is contained in a permutation in the class, then $\pi$ itself lies in the class.  Permutation classes can be specified in terms of the minimal permutations not lying in the class, which we call the {\it basis\/} of the class.  By this minimality condition, bases are necessarily {\it antichains\/}, meaning that no element of a basis is contained in another.  Although there are infinite antichains of permutations (see Atkinson, Murphy, and Ru\v{s}kuc~\cite{amr:pwocsop} for constructions and references to earlier work), we restrict our attention to finitely based classes.  Given a set of permutations $B$, we define $\Av(B)$ to be the set of permutations that avoid all of the permutations in $B$.  Thus if $\C$ is a closed class with basis $B$ then $\C=\Av(B)$, and for this reason the elements of a permutation class are often referred to as {\it restricted permutations\/}.  We let $s_n(B)$ denote the number of permutations of length $n$ in $\Av(B)$, and refer to $\sum_n s_n(B)x^n$ as the generating function of $\Av(B)$.  For more information on permutation classes, the reader is referred to B\'ona's text~\cite{bona:book}.

Each of the four systematic approaches for permutation class enumeration has a natural notion of a ``state,'' and in each case if the class is such that only finitely states are needed then --- at least in principle --- these methods give a Wilfian formula for the number of length $n$ permutations in the class.  For generating trees, the states are the labels of the isomorphic generating tree.  The classes possessing a generating tree with only finitely many labels are characterized in Vatter~\cite{finlabel}; this characterization appears here as Theorem~\ref{finlabel}.  For the insertion encoding, which associates a language to the permutation class, the natural notion of ``state'' is a state in the accepting automaton for the associated language.  The classes that require only finitely many states (or in other words, the classes that correspond to regular languages) were characterized by Albert, Linton, and Ru\v{s}kuc~\cite{insertion}; their result appears here as Theorem~\ref{insertion}.  For enumeration schemes the translation of ``state'' is ``ES$^+$-irreducible permutation'' (or, for Zeilberger's original schemes, ``ES-irreducible permutation'').  Should a class contain only finitely many such permutations then {\sc WilfPlus} can automatically enumerate it.  No characterization of these classes is known\footnote{Zeilberger~\cite{z:wilf} dismisses this by stating ``if we know beforehand that we are guaranteed to succeed, then it is not research, but doing chores.''}.  Moreover, unlike the other methods, there are subclasses of classes with finite enumeration schemes which do not themselves have finite enumeration schemes\footnote{In fact, the set of all permutations, $\Av(\emptyset)$, has a finite enumeration scheme (shown in Figure~\ref{all-perms-fig} on page~\pageref{all-perms-fig}), while several examples of classes without finite enumeration schemes are given in Section~\ref{wp-nonex}.}, indicating that such a characterization may be too much to hope for.  For substitution decompositions, simple permutations play the role of states.  As with enumeration schemes, there is no known characterization of the classes that contain only finitely many simple permutations.

\begin{figure}[t]
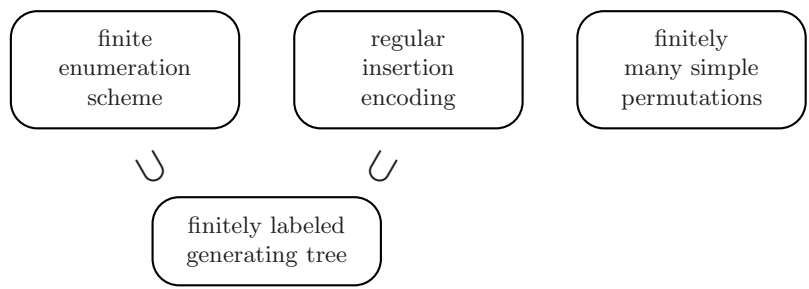

\begin{center}
\begin{psmatrix}[colsep=20pt,rowsep=10pt]
	\psframebox[linearc=10pt,cornersize=absolute]{
	\parbox{1in}{
	\begin{center}
	\begin{footnotesize}
	\begin{tabular}{c}finite\\ enumeration\\ scheme\end{tabular}
	\end{footnotesize}
	\end{center}
	}
	}
&
	\psframebox[linearc=10pt,cornersize=absolute]{
	\parbox{1in}{
	\begin{center}
	\begin{footnotesize}
	\begin{tabular}{c}regular\\ insertion\\encoding\end{tabular}
	\end{footnotesize}
	\end{center}
	}
	}
&
	\psframebox[linearc=10pt,cornersize=absolute]{
	\parbox{1in}{
	\begin{center}
	\begin{footnotesize}
	\begin{tabular}{c}finitely\\ many simple\\ permutations\end{tabular}
	\end{footnotesize}
	\end{center}
	}
	}
\\
	\multispan{2}
	\begin{psmatrix}[colsep=0pt,rowsep=10pt]
		\rput{30}{$\bigcup$}&&\rput{330}{$\bigcup$}\\
	&
		\psframebox[linearc=10pt,cornersize=absolute]{
		\parbox{1in}{
		\begin{center}
		\begin{footnotesize}
		\begin{tabular}{c}finitely labeled\\generating tree\end{tabular}
		\end{footnotesize}
		\end{center}
		}
		}
	&
	\end{psmatrix}
\end{psmatrix}
\end{center}
\caption{A depiction of the applicability of the four systematic enumeration techniques}\label{hasse-enum}
\end{figure}

The classes that these techniques can automatically enumerate are related as shown in Figure~\ref{hasse-enum}, which is to say, they are not very closely related at all (this is established via a series of examples in Section~\ref{wp-nonex} and remarks in Sections~\ref{wp-gt} and \ref{wp-ie}).  Care should be taken when reading one symbol in this diagram; while the inclusion from finitely labeled generating trees to finite enumeration schemes indicates an increase in the number of classes that can be counted, there is a corresponding decrease in information.  Finitely labeled generating trees and regular insertion encodings show that a class has a rational generating function, while classes with only finitely many simple permutations have algebraic generating functions.  It is not yet known what types of generating functions can arise from finite enumeration schemes, but they need not be algebraic.  For example, $\Av(1234)$, which has a holonomic\footnote{%
A generating function is said to be {\it holonomic\/} (or synonymously in the univariate case, {\it $D$-finite\/}) if its derivatives span a finite dimensional subspace over $\mathbb{C}(x)$.  This is equivalent to the corresponding sequence $s_n$ being holonomic (again synonymously in the univariate case, {\it $P$-recursive\/}), which means that there are polynomials $p_0,p_1,\dots p_k$ so that 
$
p_k(n)s_{n+k}+p_{k-1}(n)s_{n+k-1}+\cdots +p_0(n)s_n=0.
$
} but non-algebraic generating function (see Gessel~\cite{gessel}), has a finite enumeration scheme\footnote{%
A more trivial example would be the class of all permutations.%
} (shown in Figure~\ref{1234-fig} on page~\pageref{1234-fig}).
It is natural to hope that finite enumeration schemes produce only holonomic sequences, but this hope remains unproven.

Perhaps the greatest loss of information occurs with Wilf-equivalence.  Two classes are said to be Wilf-equivalent if they are equinumerous.  Clearly taking the reverse of a class yields a Wilf-equivalent class, as does taking the inverse, and these two operations generate the dihedral group with eight elements\footnote{%
With the exception of substitution decompositions, these techniques are not invariant under the eight permutation class symmetries.  To be precise, there are classes that cannot be handled with these methods, while their inverses can be handled easily.  Thus the comment of Albert, Linton, and Ru\v{s}kuc~\cite{insertion} that ``this apparent asymmetry does represent a possible flaw of the insertion encoding in general'' applies equally well to enumeration schemes and generating trees.%
}.
However, many examples of non-trivial Wilf-equivalences have been observed, ranging from the fact every class defined by avoiding a single pattern of length three is Wilf-equivalent\footnote{%
The classical bijective proof of this result is due to Simion and Schmidt~\cite{ss:rp}.  Zeilberger~\cite{z:snappy} gives a proof using a technique quite like enumeration schemes that generalizes to permutations of a multiset.%
}
to the theorem of Atkinson, Murphy, and Ru\v{s}kuc~\cite{amr:twostacks} that $\Av(1342)$ is Wilf-equivalent to the infinitely based class 
$$
\Av(\{2(2m-1)416385\cdots (2m)(2m-3) : m=2,3,\dots\}).
$$

If two classes both have finitely labeled generating trees, regular insertion encodings, or finitely many simple permutations then, since we can compute their generating functions from this information, we can decide whether or not they are Wilf-equivalent.  For enumeration schemes this issue is not so clear.  Occasionally, as with the enumeration schemes pictured in Figure~\ref{chow-west-fig} (a) and (b) on page~\pageref{chow-west-fig}, the Wilf-equivalence of two classes can be easily deduced from their enumeration schemes, but we present several examples in Section~\ref{wp-ex} where such deductions do not readily present themselves.

\section{Generating trees}\label{wp-gt}

Generating trees were introduced by Chung, Graham, Hoggatt, and Kleiman~\cite{cghk:baxter} and became quite popular after a pair of articles by West~\cite{west:cat, west:trees}.  The closely related {\it ECO (enumerating combinatorial objects) method\/} (see Barcucci, Del Lungo, Pergola, and Pinzani~\cite{eco:survey} for a survey) extends the notion of generating trees to other settings.

\begin{figure}[t]
\begin{footnotesize}
\begin{center}
\psset{xunit=0.03in, yunit=0.02in}
\psset{linewidth=0.25\psxunit}
\begin{pspicture}(-5,5)(155,78)
\pscircle*(10,10){1\psxunit}
\pscircle*(30,10){1\psxunit}
\pscircle*(50,10){1\psxunit}
\pscircle*(70,10){1\psxunit}
\pscircle*(90,10){1\psxunit}
\pscircle*(110,10){1\psxunit}
\pscircle*(130,10){1\psxunit}
\pscircle*(150,10){1\psxunit}
\rput[c](10,5){$1234$}
\rput[c](30,5){$4123$}
\rput[c](50,5){$3124$}
\rput[c](70,5){$4312$}
\rput[c](90,5){$2134$}
\rput[c](110,5){$4213$}
\rput[c](130,5){$3214$}
\rput[c](150,5){$4321$}
\pscircle*(20,30){1\psxunit}
\pscircle*(60,30){1\psxunit}
\pscircle*(100,30){1\psxunit}
\pscircle*(140,30){1\psxunit}
\rput[r](17,30){$123$}
\rput[l](63,30){$312$}
\rput[r](97,30){$213$}
\rput[l](143,30){$321$}
\psline(10,10)(20,30)
\psline(30,10)(20,30)
\psline(50,10)(60,30)
\psline(70,10)(60,30)
\psline(90,10)(100,30)
\psline(110,10)(100,30)
\psline(130,10)(140,30)
\psline(150,10)(140,30)
\pscircle*(40,50){1\psxunit}
\pscircle*(120,50){1\psxunit}
\rput[c](34,54){$12$}
\rput[c](126,54){$21$}
\psline(20,30)(40,50)
\psline(60,30)(40,50)
\psline(100,30)(120,50)
\psline(140,30)(120,50)
\pscircle*(80,70){1\psxunit}
\rput[c](80,76){$1$}
\psline(40,50)(80,70)
\psline(120,50)(80,70)
\end{pspicture}
\end{center}
\end{footnotesize}
\caption{The first four levels of the pattern-avoidance tree $T(132,231)$}\label{F-132-231}
\end{figure}
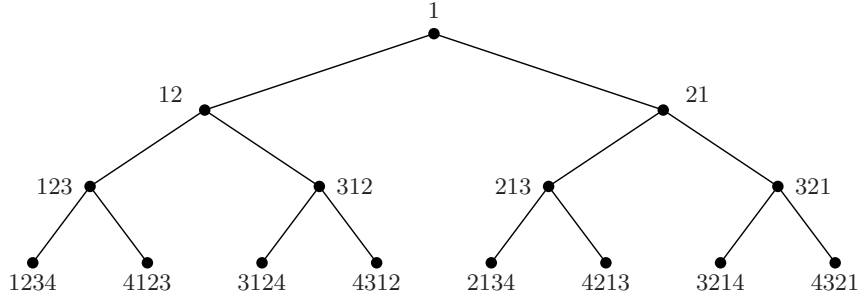

We say that the permutation $\sigma$ of length $n$ is a {\it child\/} of $\pi\in S_{n-1}$ if $\sigma$ can be obtained by inserting $n$ into $\pi$.  This defines a rooted tree $T$ on the set of all permutations.  The {\it pattern-avoidance tree\/} of $\Av(B)$, denoted by $T(B)$, is then the subtree of $T$ with nodes $\Av(B)$.  For example, the first four levels of $T(132,231)$ are shown in Figure~\ref{F-132-231}.

A generating tree, on the other hand, is a rooted, labeled tree such that the labels of the children of each node are determined by the label of that node.  Sometimes the labels of the tree are taken to be natural numbers, but this is not necessary and frequently inconvenient.  One specifies a generating tree by supplying the label of the {\it root\/} (also sometimes called the {\it axiom\/}) and a set of {\it succession rules\/} (also referred to as {\it inductive steps\/}).  For example, the complete binary tree may be given by
$$
\begin{array}{llcl}
\mbox{Root:}	&	(2)&&\\
\mbox{Rule:}	&	(2)&\leadsto &(2)(2).
\end{array}
$$

In order to enumerate the permutation class $\Av(B)$, we want to find a generating tree isomorphic (as a rooted tree) to $T(B)$.  For example, consider $T(132,231)$.  We may obtain a permutation in $\Av_n(132,231)$ by inserting $n$ either at the beginning or the end of any $\pi\in\Av_{n-1}(132,231)$, but nowhere in between, so $T(132,231)$ is isomorphic to the complete binary tree and thus to the generating tree given above.  For a more complicated example we turn to $T(1234)$, first described by West~\cite{west:cat}.  This tree is isomorphic to generating tree defined by
$$
\begin{array}{llcl}
\mbox{Root:}	&	(2,2)&&\\
\mbox{Rule:}	&	(s,t)&\leadsto &(2,t+1)(3,t+1)\cdots (s,t+1)(s,s+1)(s,s+2)\cdots (s,t)(s+1,t+t).
\end{array}
$$
While verifying this isomorphism is not difficult (consider the lexicographically first ascent and the lexicographically first occurrence of $123$ in the permutation), it is much harder to obtain the generating function for $\Av(1234)$ from this tree; for the details of this see Bousquet-M\'elou~\cite{bm:four}.

Let $T(B;\pi)$ denote the subtree of $T(B)$ that is rooted at $\pi$ and contains all descendants of $\pi$.  In an isomorphism between $T(B)$ and a generating tree, every permutation of $\Av(B)$ is assigned a label.  Clearly two permutations $\pi$ and $\sigma$ may be assigned the same label if and only if $T(B;\pi)$ and $T(B;\sigma)$ are isomorphic (again, as rooted trees).  Thus each pattern-avoidance tree $T(B)$ is isomorphic to a canonical generating tree whose labels correspond exactly to the isomorphism classes of $\{T(B;\pi) : \pi\in\Av(B)\}$.

In particular, $T(B)$ is isomorphic to a finitely labeled generating tree if and only if the set of all principal subtrees $\{T(B;\pi) : \pi\in\Av(B)\}$ contains only finitely many isomorphism classes.  When this occurs, $\Av(B)$ has a rational generating function which may be routinely computed using the transfer matrix method (see Stanley's text~\cite[Section 4.7]{stanley:ec1} for details).  The finitely based classes for which this is possible are characterized by the following theorem.

\begin{theorem}[Vatter~\cite{finlabel}]\label{finlabel}
Let $\C$ be a finitely based permutation class.  The pattern-avoidance tree of $\C$ is isomorphic to a finitely labeled generating tree if and only if $\C$ omits both a child of an increasing permutation and a child of a decreasing permutation.
\end{theorem}

For example, $T(132,231)$ satisfies the hypotheses of Theorem~\ref{finlabel} because it omits both $132$ (a child of the increasing permutation $12$) and $231$ (a child of $21$).  A less trivial example is given by $T(123,3214,2143,15432)$, which arose in Klazar~\cite{k:growth}.  The Maple package {\sc FinLabel} (described in \cite{finlabel} and available at \url{http://math.rutgers.edu/~vatter/}) can find the generating functions for classes satisfying Theorem~\ref{finlabel} completely automatically.

It is easy to see that the hypotheses of Theorem~\ref{finlabel} are necessary%
\footnote{Suppose, without loss, that $\C$ contains all children of every increasing permutation.  Then for each $n$, $12\cdots n$ has $n+1$ children in the pattern-avoidance tree of $\C$, and thus no two of these nodes may share the same label.}%
.  The other direction is proved by showing that every sufficiently long permutation is ``GT-reducible.''  Since GT-reducibility is a stronger condition than the ES-reducibility of enumeration schemes, every class with a finitely labeled generating tree has a finite enumeration scheme.

\section{The insertion encoding}\label{wp-ie}

The insertion encoding, recently introduced by Albert, Linton, and Ru\v{s}kuc~\cite{insertion}, is a correspondence between permutation classes and languages.  With it, one may attack the enumeration problem with all the tools of formal language theory.  Roughly, this correspondence associates to each permutation a word describing how that permutation evolved.  At each stage until the desired permutation has been constructed, at least one open {\it slot\/} (represented by a $\diamond$) exists in the intermediate {\it configuration\/}, and to proceed to the next configuration we insert a new maximal entry into one of these slots.  This insertion can occur in four possible ways:
\begin{itemize}
\item the slot can be filled (replacing a $\diamond$ by $n$),
\item the new entry can be inserted to the left of the slot (replacing a $\diamond$ by $n\,\diamond$), 
\item the new entry can be inserted to the right (replacing a $\diamond$ by $\diamond\,n$), or
\item the slot can be divided into two slots with the new entry in between (replacing a $\diamond$ by $\diamond\,n\,\diamond$).
\end{itemize}
These operations are denoted by the symbols $\ie{f}$, $\ie{l}$, $\ie{r}$, and $\ie{m}$, respectively.  Since each of these operations can be performed on any open slot at any stage, we subscript their symbols with the number of the slot they were applied to (read from left to right).  For example, the permutation $31254$ has the insertion encoding $\ie{m}_1\ie{l}_2\ie{f}_1\ie{r}_1\ie{f}_1$ because its evolution is
$$
\begin{array}{c}
\diamond\\
\diamond\,1\,\diamond\\
\diamond\,12\,\diamond\\
312\,\diamond\\
312\,\diamond\,4\\
31254
\end{array}
$$

Let $\sb(k)$ denote the permutation class whose basis consists of all length $2k+1$ permutations of the form $babab\cdots bab$ where the $a$'s represent the elements $\{1,2,\dots,k\}$ and the $b$'s represent the elements $\{k+1,k+2,\dots,2k+1\}$.  These classes are called {\it slot bounded\/} because in the evolution of a permutation in $\sb(k)$ there are never more than $k$ open slots.

\begin{theorem}[Albert, Linton, and Ru\v{s}kuc~\cite{insertion}]\label{insertion}
The insertion encoding of a finitely based class is regular if and only if the class is a subclass of $\sb(k)$ for some $k$.
\end{theorem}

One can show using the Erd\H{o}s-Szekeres theorem~\cite{es:acpig} (or one can refer to the proof in \cite{insertion}) that Theorem~\ref{insertion} includes all of the classes identified by Theorem~\ref{finlabel} as having finitely labeled generating trees.

Even when the insertion encoding of a class is not regular, useful information can still be obtained by this correspondence.  For example, Albert, Elder, Rechnitzer, Westcott, Zabrocki~\cite{1324} used regular approximations to the insertion encoding of $\Av(1324)$ to establish that $s_n(1324)>9.35^n$ for sufficiently large $n$, thereby disproving a conjecture of Arratia~\cite{arratia}.  Additionally, Albert, Linton, and Ru\v{s}kuc~\cite{insertion} consider several classes with context-free insertion encodings and are able to obtain their (algebraic) generating functions from these languages.  However, the derivation of insertion encodings is only automatic for subclasses of $\sb(k)$, and thus we choose to limit our focus to this case.

\section{Substitution decompositions}\label{wp-simple}

Substitution decompositions (also known as modular decompositions, disjunctive decompositions, and $X$-joins) have proven to be a useful technique in a wide range of settings, ranging from game theory to combinatorial optimization (see M\"ohring~\cite{m:aas} or M\"ohring and Radermacher~\cite{mr:sdd} for extensive references).  Permutation class enumeration is no exception.

An {\it interval\/} (also called a {\it block\/}, or in other contexts, {\it factor\/}, {\it clan\/}, or even {\it convex subset\/}) in the permutation $\pi$ is an interval of indices $I=[a,b]$ such that the set of values $\{\pi(i) : i\in I\}$ also forms an interval.  Clearly every permutation of length $n$ has $n$ trivial intervals of length one and one trivial interval of length $n$.  A permutation that has no non-trivial intervals is called {\it simple\/} (the analogous term in other contexts is often {\it prime\/} or {\it primitive\/}).  

Simple permutations first appear in the work of Atkinson and Stitt~\cite{as:wreath}, which is followed up by Albert and Atkinson~\cite{aa:simple}.  Although in other contexts substitution decompositions are most often applied to algorithmic problems, they also have powerful enumerative applications.  A class with only finitely many simple permutations has a recursive structure in which long permutations are built up from smaller permutations (their intervals).  Thus it is natural to expect these classes to have algebraic generating functions, and this intuition is borne out by the following theorem.

\begin{theorem}[Albert and Atkinson~\cite{aa:simple}]\label{simple}
A permutation class with only finitely many simple permutations has an algebraic generating function.
\end{theorem}

The canonical example of a class with only finitely many simple permutations is $\Av(132)$.  By considering the entries to the left and to the right of the $n$ in a permutation in $\Av_n(132)$ one simultaneously derives a decomposition of these permutations that leads immediately to the Catalan numbers and sees that this class contains no simple permutations of length three\footnote{Actually, there are no simple permutations of length three, $132$-avoiding or otherwise.} or longer.

Another example of a class with only finitely many simple permutations is the class of {\it separable permutations\/}.  This class, first introduced by Bose, Buss, and Lubiw~\cite{bose:matching}, is essentially the permutation analogue of series-parallel posets (see Stanley~\cite{s:epg,stanley:ec1}) and complement reducible graphs (see Corneil, Lerchs, and Burlingham~\cite{clb:crg}).  To define separable permutations we first need two binary operations on permutations.  Given two permutations $\pi\in S_m$ and $\sigma\in S_n$ we define their {\it direct sum\/}, written $\pi\oplus\sigma$, by
$$
(\pi\oplus\sigma)(i)
=
\left\{
\begin{array}{ll}
\pi(i)&\mbox{if $i\in [m]$,}\\
\sigma(i-m)+m&\mbox{if $i\in[m+n]\setminus[m]$.}
\end{array}
\right.
$$
Similarly, we define their {\it skew sum\/}, $\pi\ominus\sigma$, by
$$
(\pi\ominus\sigma)(i)
=
\left\{
\begin{array}{ll}
\pi(i)+n&\mbox{if $i\in [m]$,}\\
\sigma(i-n)&\mbox{if $i\in[m+n]\setminus[m]$.}
\end{array}
\right.
$$
Given a class $\C$, we denote by $\strongcomp(\C)$ the {\it strong completion\/} of $\C$, which is the smallest class containing $\C$ such that both $\pi\oplus\sigma$ and $\pi\ominus\sigma$ lie in $\strongcomp(\C)$ for every $\pi,\sigma\in\strongcomp(\C)$.

The separable permutations are the strong completion of $\{1\}$.  As was shown by Bose, Buss, and Lubiw~\cite{bose:matching}, this class can also be described as $\Av(2413,3142)$.  The enumeration of this class (which can now be seen to follow routinely from Theorem~\ref{simple} and the fact that the only simple separable permutations are $1$, $12$, and $21$) was first undertaken by West~\cite{west:cat}.  He used generating trees to show that the separable permutations are counted by the large Schr\"oder numbers.  Later, Ehrenfeucht, Harju, ten Pas, and Rozenberg~\cite{ehpr:schroeder} (who also gave another proof that the basis of this class is $\{2413,3142\}$) presented a bijection between separable permutations and parenthesis words, the objects Schr\"oder was originally interested in counting.

One of the notable features of Theorem~\ref{simple} is that it does not seem to require the class to be finitely based.  However, this is merely an illusion:

\begin{theorem}[Albert and Atkinson~\cite{aa:simple}, Murphy~\cite{maximillian}]\label{simplepwo}
A permutation class with only finitely many simple permutations is both finitely based and partially well-ordered\,\footnote{A partially ordered set is said to be partially well-ordered if contains neither an infinite strictly decreasing subsequence (which is never possible for a permutation class) nor an infinite antichain.}.
\end{theorem}

There is a semi-algorithm for establishing that a class contains only finitely many simple permutations.  This semi-algorithm stems from the following theorem of Schmerl and Trotter~\cite{st:simple}, who proved it in the more general context of binary relational systems.  Versions of the theorem for $2$-structures and $k$-structures are given by Ehrenfeucht and Rozenberg~\cite{er:ph2s} and Ehrenfeucht and McConnell~\cite{em:kgt}, repectively, and a proof for the special case of permutations can be found in Murphy's thesis~\cite{maximillian}.

\begin{theorem}[Schmerl and Trotter~\cite{st:simple}]\label{simple:hered}
Every simple permutation of length $n>2$ contains a simple permutation of length $n-1$ or $n-2$.
\end{theorem}

If the class $\C$ contains only finitely many simple permutations, then clearly there is an integer $n$ so that $\C$ does not contain any simple permutations of lengths $n-1$ or $n-2$.  In the other direction, Theorem~\ref{simple:hered} shows that if we have found such an integer $n$ then $\C$ contains no simple permutations of length $n-2$ or longer.  Therefore, when a class happens to contain only finitely many simple permutations, this fact can be verified automatically.

It remains an interesting open question if it is decidable whether a class contains only finitely many simple permutations.

\section{Zeilberger's original enumeration schemes}\label{wp-wilf}

Zeilberger~\cite{z:wilf} developed the notion of {\it enumeration schemes\/} and wrote the Maple package {\sc Wilf} to automate their discovery.  Roughly, enumeration schemes are a divide and conquer technique which aims to partition the class into smaller pieces from which recurrences can be derived.

Take $\pi\in S_k$, suppose that $n\ge k$, and let $1\le i_1<i_2<\cdots<i_k\le n$.  In Zeilberger's original formalization of enumeration schemes, we divide $\Av_n(B)$ into the sets
$$
A_\pi(n;B;i_1,i_2,\dots, i_k)=\{p\in\Av_n(B) : p(1)=i_{\pi(1)},\dots,p(k)=i_{\pi(k)}\}.
$$
In words, $A_\pi(n;B;i_1,i_2,\dots, i_k)$ is the set of $B$-avoiding length $n$ permutations that begin with the entries $i_1,i_2,\dots,i_k$, in the order specified by $\pi$.  For example,
\begin{equation}
A_{312}(9;B;2,3,7)=\{723x_4x_5x_6x_7x_8x_9\in\Av_n(B)\}.
\label{z-set-example}
\end{equation}

In order to make enumeration schemes more closely resemble generating trees and the insertion encoding, we consider a symmetry of his approach.  Everywhere Zeilberger mentions a permutation we consider its inverse.  Thus we should specify the set of restrictions, $B$, a set of small entries of some length, $\pi$, and the positions in which the entries of $\pi$ occur.  But instead of specifying the positions, we specify the gaps between the entries with a {\it gap vector\/}, $\g$.  After performing these transformations, our version of \eqref{z-set-example} is
$$
Z(B;231;(1,0,3,2))=\{x_123x_4x_5x_61x_8x_9\in\Av_9(B)\},
$$
and in general we are concerned with the sets
$$
Z(B;\pi;\g)=\{p\in\Av_{k+\|\g\|}(B) : p(g_1+1)=\pi(1),\dots,p(g_1+\cdots+g_{k}+k)=\pi(k)\},
$$
where $k$ is the length of $\pi$ and $\|\g\|$ denotes the sum of the components of $\g$.
Thus $Z(B;\pi;\g)$ is the set of all $B$-avoiding permutations of length $k+\|\g\|$ whose least $k$ elements occur in the positions $g_1+1,g_1+g_2+2,\dots,g_1+g_2+\dots+g_k+k$ and form a $\pi$-subsequence.

Not all pairs $(\pi,\g)$ result in a nonempty $Z$-set.  Following Zeilberger, for a length $k$ permutation $\pi$ we define
$$
\J(\pi)=\{j\in[k+1] : Z(B;\pi;\g)=\emptyset\mbox{ for all $\g$ with $g_j>0$}\}.
$$
Thus $Z(B;\pi;\g)$ is guaranteed to be empty if $\g$ does not ``obey'' $\J(\pi)$, meaning that $g_j\neq 0$ for some $j\in\J(\pi)$.

For example, consider the case $B=\{132\}$.  Then $2\in \J(12)$ because if $g_2>0$ then there is some entry between $1$ and $2$ in every permutation in $Z(B;\pi;\g)$, and this gives a $132$-pattern.  In order to check that $\J(12)=\{2\}$ we need merely observe that $312$ and $123$ avoid $132$.  Our approach in this example can easily be generalized to compute $\J(\pi)$ for any $\pi$ and $B$.

\begin{proposition}\label{J-computable}
For any permutation $\pi$ and basis $B$, $\J(\pi)$ can be computed by inspecting the $B$-avoiding children of $\pi$.
\end{proposition}
\begin{proof}
Consider the vector $\h$ for which $h_i=0$ for all $i\neq j$ and $h_j=1$.  If $Z(B;\pi;\h)=\emptyset$  then $Z(B;\pi;\g)=\emptyset$ for all $\g$ with $g_j>0$, so $j\in\J(\pi)$.  If instead $Z(B;\pi;h)\neq\emptyset$ then $j\notin \J(\pi)$.
\end{proof}

For any $r\in[k]$, the set $Z(B;\pi;(g_1,\dots,g_{k+1}))$ embeds naturally (remove the entry $\pi(r)$ and standardize) into
\begin{equation}\label{embedding}
Z(B;\st(\pi-\pi(r));(g_1,\dots,g_{r-1},g_r+g_{r+1},g_{r+2},\dots,g_{k+1})),
\end{equation}
where $\pi-\pi(r)$ denotes the word obtained from $\pi$ by omitting the entry $\pi(r)$, so, for example, $51342-1=5342$.  To make \eqref{embedding} easier to state, we define $\d_r(\pi)$ to be $\st(\pi-\pi(r))$ and let
$$
\d_r((g_1,\dots,g_{k+1}))=(g_1,\dots,g_{r-1},g_r+g_{r+1},g_{r+2},\dots,g_{k+1}).
$$

Sometimes the embedding of $Z(B;\pi;\g)$ into $Z(B;\d_r(\pi);\d_r(\g))$ is a bijection.  If this is true for all gap vectors $\g$ that obey $\J(\pi)$, that is, that have $g_j=0$ for all $j\in\J(\pi)$, then we say that $\pi(r)$ is {\it enumeration-scheme-reducible for $\pi$ with respect to $B$\/}, or, for short, ES-reducible.  (Zeilberger~\cite{z:wilf} refers to such entries as {\it reversely deleteable\/}.)  We also say that a permutation with an ES-reducible entry is itself ES-reducible, and a permutation without an ES-reducible entry is ES-irreducible.

For example, suppose again that $B=\{132\}$ and consider the permutation $12$.  We have already observed that $\J(12)=\{2\}$.  Now we claim that the entry $1$ is ES-reducible.  The gap vectors that obey $\J(12)$ are those of the form $(g_1,0,g_3)$, and thus we would like to verify that the embedding of $Z(\{132\};12;(g_1,0,g_3))$ into $Z(\{132\};1;(g_1,g_3))$ is a bijection.

Take $p\in Z(\{132\};1;(g_1,g_3))$ and consider inverting this embedding.  In this case, that amounts to inserting the element $1$ into position $g_1+1$ and increasing all other entries of $p$ by $1$.  Label the resulting permutation $p'$.  For example, from $\g=(3,0,1)$ and $p=52314$ we obtain $p'=634125$.

We would like to show that $p'$ avoids $132$.  To show this we consider all possible ways in which the new element $1$ could participate in a $132$-pattern.  Clearly this entry must be the first entry in such a pattern.  Now note that the $2$ in $p'$ cannot participate in this $132$-pattern, because the $1$ and $2$ are adjacent.  But then there is a $132$-pattern in $p'$ which uses the $2$ instead of the $1$, and thus $p$ contains a $132$-pattern, a contradiction.

Thus we have shown that
$$
|Z_n(\{132\};12;(g_1,g_2,g_3))|
=
\case{0}{if $g_2>0$,}
{|Z_{n-1}(\{132\};1;(g_1,g_3))|}{if $g_2=0$.}
$$

Although this example did not demonstrate it, detecting and verifying ES-reducibility by hand can be enormously tedious.  Fortunately, it is also unnecessary.  By adapting the approach used in \cite{finlabel}, we arrive at the following test for ES-reducibility that can be routinely checked by computer%
\footnote{Zeilberger's approach in \cite{z:wilf} used what he referred to as ``logical reasoning,'' and while it is no less rigorous than this approach, Proposition~\ref{test-rd} has the advantage of being very explicit.}%
. In it we let $\Binf$ denote the length of the longest permutation in $B$.

\begin{proposition}\label{test-rd}
The entry $\pi(r)$ of the permutation $\pi$ is ES-reducible if and only if
$$
|Z(B;\pi;\g)|=|Z(B;\d_r(\pi);\d_r(\g))|
$$
for all gap vectors $\g$ of the appropriate length that obey $J(\pi)$ and satisfy $\|\g\|\le\Binf-1$.
\end{proposition}
\begin{proof}
If $\pi(r)$ is ES-reducible then the claim follows by definition.  To establish the other direction, suppose that $\pi(r)$ is not ES-reducible, and choose $\g$ and $p\in Z(B;\d_r(\pi);\d_r(\g))$ so that $\g$ obeys $\J(\pi)$ but $p$ cannot be obtained from a permutation in $Z(B;\pi;\g)$ by removing $\pi(r)$ and standardizing.

First form the ($B$-containing) permutation $p'$ by incrementing each entry of $p$ that is at least $\pi(r)$ by $1$ and inserting $\pi(r)$ into position $g_1+\cdots+g_r+r$.  Thus $p'$ is the permutation that would have mapped to $p$, except that $p'$ contains a pattern from $B$ and thus does not lie in $Z(B;\pi;\g)$.

Now pick some $\beta\in B$ that is contained in $p'$, and choose a specific occurrence of $\beta$ in $p'$.  Note that since $p=\st(p'-\pi(r))$ avoids $B$, this occurrence of $\beta$ must include the entry $\pi(r)$.  Let $p''$ denote the standardization of the subsequence of $p'$ formed by all entries that are either in the chosen occurrence of $\beta$ or in $\pi$ (or in both), so $p''$ contains a $\beta$-pattern and lies in $Z(\emptyset;\pi;\h)$ for some $\J(\pi)$-obeying $\h$ with $\|\h\|\le\Binf-1$.  On the other hand, $\st(p''-\pi(r))$ avoids $B$, which implies that $|Z(B;\d_r(\pi);\h)|>|Z(B;\pi;\h)|$, as desired.
\end{proof}

For example, consider the basis $B=\{132\}$ again.  In order to show that the $1$ in $\pi=12$ is ES-reducible using this proposition, we first find that $\J(12)=\{2\}$ and then perform the following 10 computations:

\begin{footnotesize}
$$
\begin{array}{ccc}
\g&|Z(\{132\};12;\g)|&|Z(\{132\};1;d_1(\g))|\\\hline
(0,0,0)&1&1\\
(0,0,1)&1&1\\
(1,0,0)&1&1\\
(0,0,2)&1&1\\
(1,0,1)&2&2\\
(2,0,0)&2&2\\
(0,0,3)&1&1\\
(1,0,2)&3&3\\
(2,0,1)&5&5\\
(3,0,0)&5&5
\end{array}
$$
\end{footnotesize}

In a similar manner one can verify that the $1$ in $21$ is ES-reducible and that $\J(21)=\emptyset$.  This gives the following enumeration scheme for $\Av(132)$:

\begin{eqnarray*}
s_n(132)
&=&
|Z(\{132\}; \emptyset; (n))|,
\\
|Z(\{132\}; \emptyset; (g_1))|
&=&
\sum_{i=0}^{g_1-1} |Z(\{132\}; 1; (i, g_1-i-1))|,
\\
|Z(\{132\}; 1; (g_1,g_2))|
&=&
\sum_{i=0}^{g_1-1} |Z(\{132\}; 21; (i, g_1-i-1, g_2))|
\\
&&
\quad+\sum_{i=0}^{g_2-1} |Z(\{132\}; 12; (g_1,i,g_2-i-1))|,
\\
&=&
\sum_{i=0}^{g_1} |Z(\{132\}; 1; (i, g_1+g_2-i-1))|.
\end{eqnarray*}

\section{Extending enumeration schemes}\label{wp-wp}

We will replace the sets $\J(\pi)$ in this section, giving us a more powerful version of enumeration schemes that can be found automatically with the Maple package {\sc WilfPlus}.

In order to motivate this change, we first consider a shortcoming of $\J(\pi)$.  Let $B=\{1342,1432\}$.  It can be shown easily, even by hand, that $12$ is ES-irreducible.  To do so, first note that $\J(12)=\emptyset$, as witnessed by the permutations $312$, $132$, and $123$.  Now consider the set $Z(B;12;(0,2,0))$.  This set is empty, but removing $1$ gives the nonempty set $Z(B;1;(2,0))$ while removing $2$ gives the nonempty set $Z(B;1;(0,2))$.  Indeed, this reasoning generalizes to show that all permutations of the form $\ominus^m 12$ are ES-irreducible, so $\Av(1342,1432)$ does not have a finite enumeration scheme, at least in Zeilberger's original sense.

Zeilberger's enumeration schemes fail in the previous example for a very simple reason: $\J(12)$ is too coarse to capture the fact that $(0,2,0)$ is not a valid gap vector for a $B$-avoiding descendant of $12$.  We remedy this problem with the following definition.

\begin{definition}
The entry $\pi(r)$ of the length $k$ permutation $\pi$ is said to be {\it \esp-reducible} if
\begin{eqnarray*}\label{wrd-def}
|Z(B;\pi;\g)|=|Z(B;\d_r(\pi);\d_r(\g))|
\end{eqnarray*}
whenever $Z_n(B;\pi;\g)$ is nonempty.  Further, we say that the permutation $\pi$ is \esp-reducible if it contains an \esp-reducible entry, and \esp-irreducible otherwise.
\end{definition}

We then replace (for now) the set $\J$ by
$$
\G(\pi) = \{\g : Z(B;\pi;\g)\neq\emptyset\}.
$$

Recall that for $\pi(r)$ to be ES-reducible, it had to satisfy (\ref{wrd-def}) for all $\g$ that contained $0$'s in the positions specified by $\J(\pi)$.  Clearly if $\g\in\G(\pi)$ then $\g$ obeys $\J(\pi)$, but there can be gap vectors that obey $\J(\pi)$ and do not lie in $\G(\pi)$, as in our previous example with $B=\{1342,1432\}$.  Thus we have obtained a weaker condition by requiring the satisfaction of (\ref{wrd-def}) less often.

The proof of Proposition~\ref{test-rd} carries over to this context to give the analogous result on testing for \esp-reducibility.

\begin{proposition}\label{test-wrd}
The entry $\pi(r)$ of the permutation $\pi$ is \esp-reducible if and only if
$$
|Z(B;\pi;\g)|=|Z(B;\d_r(\pi);\d_r(\g))|
$$
for all $\g\in\G(\pi)$ with $\|\g\|\le\Binf-1$.
\end{proposition}

One can view $\G(\pi)$ as an lower order ideal\footnote{This means that $\x\in\G(\pi)$ whenever $\x\le \y$ for some $\y\in\G(\pi)$.} of $\mathbb{N}^{|\pi|+1}$ under the product order, where $(x_1,\dots,x_k)\le (y_1,\dots,y_k)$ if and only if $x_i\le y_i$ for all $i\in[k]$.  Therefore we carry our definitions about permutation classes over to this context.  In particular, we say that the {\it basis\/} of $\G(\pi)$ is the set of minimal vectors not in $\G(\pi)$, and if $B$ is a set of vectors then we write $\Av(B)$ to denote the set $\{\g : \g\not\ge\b\mbox{ for all }\b\in B\}$.

For example, let us compute the basis of $\G(12)$ when $B=\{1342,1432\}$.  As already observed, $(0,2,0)$ does not lie in $\G(12)$.  This gap vector is minimal in $\mathbb{N}^3\setminus \G(12)$ because $Z(B;12;(0,1,0))$ is nonempty.  To show that the basis of $\G(12)$ is precisely $(0,2,0)$, it suffices to note that the permutation
$$
3\ 4\ \cdots\ (g_1+2)\ 1\ (g_1+3)\ 2\ (g_1+4)\ (g_1+5)\ \cdots\ (g_1+g_2+3)
$$
avoids $B$, so $(g_1,1,g_2)\in\G(12)$ for all $g_1,g_2\in\mathbb{N}$.  Thus we have shown that $\G(12)=\Av((0,2,0))$.

Now that we have computed $\G(12)$, it is not hard to check that $2$ is \esp-reducible for $12$.  In order to do so we need to show that the embedding in question is a bijection for all gap vectors $(g_1,g_2,g_3)$ with $g_2\le 1$.

Suppose to the contrary that the embedding is not onto and take $p\in Z(B;1;(g_1,g_2+g_3))$ that is not mapped to.  In other words, the permutation $p'$ obtained from $p$ by inserting $2$ into position $g_1+g_2+2$ and incrementing all the entries of $p$ of value at least $2$ contains a $1342$ or $1432$-pattern.  Since $p$ avoids $\{1342,1432\}$, this pattern must involve the entry $2$.  First, the $2$ cannot play the role of the ``$2$'' in such a pattern, because then the $1$ would be forced to play the role of the ``$1$,'' and there can be at most one entry between the $1$ and $2$ since $g_2\le 1$.  The only other possible role for the $2$ is as the ``$1$,'' but in this case we could substitute the $1$, thereby finding a $B$-pattern in $p$, a contradiction.

Of course, we can replace $\J(\pi)$ by $\G(\pi)$ only if we are able to work with $\G(\pi)$.  Because $\mathbb{N}^{|\pi|+1}$ is partially well-ordered by the product ordering\footnote{This fact, which is not difficult to prove, can be found many places, for example Nash-Williams~\cite{nw:wqo}.}, the basis of $\G(\pi)$, which is by definition an antichain, must be finite.  However, this basis may be quite large, or may contain vectors with large components.  For example, another way to state the Erd\H{o}s-Szekeres theorem~\cite{es:acpig} is that
$$
\G(\emptyset)=\Av(((j-1)(k-1)+1))
$$
when $B=\{12\cdots j, k\cdots 21\}$.  While this does not preclude effective computation of $\G(\pi)$, it does suggest that such computations could be time consuming.

It happens that we can circumvent this problem by replacing $\G(\pi)$ by a different set of gap vectors.  First we return to the way in which $\G(\pi)$ is used.  By our definitions, if $\pi(r)$ is \esp-reducible for $\pi$ then
\begin{eqnarray*}
|Z(B;\pi;\g)|=
\case{0}{if $\g\notin\G(\pi)$,}
{|Z(B;\d_r(\pi);\d_r(\g))|}{otherwise.}
\end{eqnarray*}
This equality shows that when enumerating the $B$-avoiding descendants of $\pi$ with gap vector $\g$, we first check to see if $\g$ lies in $\G(\pi)$.  If $\g\notin\G(\pi)$ then we can be sure that no such descendants exist.  If $\g\in\G(\pi)$, then $Z(B;\pi;\g)$ is in one-to-one correspondence with $Z(B;\d_r(\pi);\d_r(\g))$.  Note that if $Z(B;\pi;\g)$ and $Z(B;\d_r(\pi);\d_r(\g))$ are both empty then they are trivially in one-to-one correspondence, so we could instead use the recurrence
\begin{eqnarray*}
|Z(B;\pi;\g)|=
\case{0}{if $\g\notin\G(\pi)$ and $\d_r(\g)\in\G(\d_r(\pi))$,}
{|Z(B;\d_r(\pi);\d_r(\g))|}{otherwise.}
\end{eqnarray*}

This equality shows that instead of considering $\G(\pi)$, we can look at the larger set of gap vectors for which either $Z(B;\pi;\g)\neq\emptyset$ or $Z(B;d_r(\pi);d_r(\g))=\emptyset$.  Unfortunately, this set need not be an ideal%
\footnote{An example of this occurs with $B=\{231,k\cdots 21\}$.  In Proposition~\ref{prop-chow-west} we observe that the set of $\g\in\mathbb{N}^3$ for which $Z(B;21;\g)\neq\emptyset$ is $\Av((0,1,0),(k-2,0,0))$ while the set of $\g\in\mathbb{N}^2$ for which $Z(B;1;\g)=\emptyset$ is $\Av((k-1,0))$, so the set of $\g\in\mathbb{N}^3$ for which either $Z(B;21;\g)\neq\emptyset$ or $Z(B;1;\d_1(\g))=\emptyset$ is $\Av((0,1,0))\setminus\{(k-2,0,0)\}$.},
so we consider the largest lower order ideal of $\mathbb{N}^{|\pi|+1}$ for which these conditions hold:
$$
\G_r(\pi)=\{\g\in\mathbb{N}^{|\pi|+1} : Z(B;\pi;\h)\neq\emptyset\mbox{ or }Z(B;d_r(\pi);d_r(\h))=\emptyset\mbox{ for all }\h\le\g\}.
$$
Note that $\G(\pi)\subseteq\G_r(\pi)$: if $\g\in\G(\pi)$ then $Z(B;\pi;\g)\neq\emptyset$, so $Z(B;\pi;\h)\neq\emptyset$ for all $\h\le\g$, so $\g\in\G_r(\pi)$.  With this observation we have
\begin{eqnarray*}
|Z(B;\pi;\g)|=
\case{0}{if $\g\notin\G_r(\pi)$,}
{|Z(B;\d_r(\pi);\d_r(\g))|}{otherwise,}
\end{eqnarray*}
if $\pi(r)$ is \esp-reducible for $\pi$ and $B$.

This new set has several advantages over $\G(\pi)$.  For one, it may be considerably smaller, thus simplifying the scheme.  More importantly, the following bound on basis elements implies that $\G_r(\pi)$ can be found automatically.

\begin{proposition}\label{Hbasis}
Suppose that $\pi(r)$ is an \esp-reducible entry in $\pi$ (with respect to the permutation class with basis $B$).  Then each basis vector $\b$ of $\G_r(\pi)$ satisfies $\|\b\|\le\Binf-1$.
\end{proposition}
\begin{proof}
Suppose to the contrary that $\G_r(\pi)$ has a basis vector $\b$ with $\|\b\|\ge\Binf$.  Then $Z(B;\pi;\a)=\emptyset$ and $Z(B;\d_r(\pi);\d_r(\a))\neq\emptyset$ for some $\a\le\b$ since $\b\notin \G_r(\pi)$.  However, because $\b$ is a basis vector for $\G_r(\pi)$, every $\a<\b$ lies in $\G_r(\pi)$, so 
$Z(B;\pi;\a)\neq\emptyset$ or $Z(B;\d_r(\pi);\d_r(\a))=\emptyset$ for these vectors.  If $Z(B;\d_r(\pi);\d_r(\a))=\emptyset$ for one of these vectors then $Z(B;\d_r(\pi);\d_r(\b))$ must also be empty, but this contradicts the fact that $\b\notin\G_r(\pi)$, so $Z(B;\pi;\a)\neq\emptyset$ for all $\a<\b$.  Therefore $Z(B;\pi;\b)=\emptyset$ and $Z(B;\d_r(\pi);\d_r(\b))\neq\emptyset$.

Now, since $\pi(r)$ is \esp-reducible and $Z(B;\pi;\a)\neq\emptyset$ for all $\a<\b$, we know that there is a bijection between $Z(B;\pi;\a)$ and $Z(B;\d_r(\pi);\d_r(\a))$ for all $\a<\b$.  Thus we have the following diagram.

$$
\begin{psmatrix}[colsep=1.5cm,rowsep=1.5cm]
Z(B;\pi;\b)=\emptyset&Z(B;d_r(\pi);d_r(\b))\neq\emptyset\\
Z(B;\pi;\a)\neq\emptyset&Z(B;d_r(\pi);d_r(\a))\neq\emptyset
\psset{arrows=->,labelsep=3pt,nodesep=3pt}
\ncline{1,1}{2,1}
\ncline{1,2}{2,2}
\psset{arrows=<->}
\ncline{2,1}{2,2}
\end{psmatrix}
$$

The rest of the proof is similar to the proof of Proposition~\ref{test-rd}.  Choose a permutation $p\in Z(B;\d_r(\pi);\d_r(\b))$ and form $p'$ by incrementing each entry of $p$ that is at least $\pi(r)$ by $1$ and inserting $\pi(r)$ into position $g_1+\cdots+g_r+r$.  Choose a specific occurrence of some $\beta\in B$ in $p'$.  Since $p$ avoids $B$, this occurrence of $\beta$ must involve the entry $\pi(r)$.  Let $p''$ denote the standardization of the subsequence of $p'$ given by the entries from $\pi$ together with the entries from the chosen occurrence of $\beta$.  Therefore $p''$ contains a $\beta$-pattern and lies in $Z(\emptyset;\pi;\a)$ for some $\a<\b$ with $\|\a\|\le\Binf-1$.  However, this is a contradiction because $\d_r(p'')$ avoids $B$ and thus lies in $Z(B;\d_r(\pi);\d_r(\a))$, and we have assumed that $\d_r$ is a bijection between $Z(B;\pi;\a)$ and $Z(B;\d_r(\pi);\d_r(\a))$.
\end{proof}

We conclude this section by writing out the enumeration scheme that we have derived for $\Av(1342,1432)$.

\begin{eqnarray*}
s_n(1342,1432)
&=&
|Z(\{1342,1432\}; \emptyset; (n))|,
\\
|Z(\{1342,1432\}; \emptyset; (g_1))|
&=&
\sum_{i=0}^{g_1-1} |Z(\{1342,1432\}; 1; (i, g_1-i-1))|,
\\
|Z(\{1342,1432\}; 1; (g_1,g_2))|
&=&
\sum_{i=0}^{g_1-1} |Z(\{1342,1432\}; 21; (i, g_1-i-1, g_2))|
\\
&&
\quad+\sum_{i=0}^{g_2-1} |Z(\{1342,1432\}; 12; (g_1,i,g_2-i-1))|,
\\
&=&
\sum_{i=0}^{g_1-1} |Z(\{1342,1432\}; 1; (i, g_1+g_2-i-1))|
\\
&&
\quad+2|Z(\{1342,1432\}; 1; (g_1,g_2-1))|.
\end{eqnarray*}

\section{A collection of failures}\label{wp-nonex}

In this section we collect numerous negative results needed to justify the lack of inclusions in Figure~\ref{hasse-enum}.  The class $\Av(123)$ shows immediately that classes can have finite enumeration schemes without having regular insertion encodings.  An example of a class with a regular insertion encoding but without a finite enumeration scheme is given by $\Av(1234,4231)$. It can be computed that this classes lies in $\sb(4)$, but the following proposition shows that it does not have a finite enumeration scheme.

\begin{proposition}\label{1234-4231-bad}
For all $k$, the permutation $k\cdots 21$ is \esp-irreducible for $\Av(1234,4231)$.
\end{proposition}
\begin{proof}
Let $\pi=k\cdots 21$.  First we show that the entries $\pi(r)=k-r+1$ for $r\in[k-1]$ are \esp-irreducible.  Let $\g$ denote the vector in $\mathbb{N}^{k+1}$ which is identically zero except for $g_{k-r+1}$ and $g_{k-r+2}$, which are both $1$ (these two components correspond to the gaps on either side of $\pi(k-r+1)$).  Now observe that there are two permutations in $Z(\{1234,4231\};\d_r(\pi);\d_r(\g))$:
$$
\begin{array}{ll}
\st(k(k-1)\cdots (k-r+2)(k+1)(k+2)(k-r)\cdots 21),&\mbox{and}\\
\st(k(k-1)\cdots (k-r+2)(k+2)(k+1)(k-r)\cdots 21),
\end{array}
$$
while $Z(\{1234,4231\};\pi;\g)$ contains only one permutation,
$$
k(k-1)\cdots (k-r+2)(k+1)(k-r+1)(k+2)(k-r)\cdots 21.
$$
Thus $\g\in\G_r(\pi)$ but $Z(\{1234,4231\};\pi;\g)$ and $Z(\{1234,4231\};\d_r(\pi);\d_r(\g))$ are not in one-to-one correspondence, so $\pi(r)$ is not \esp-reducible for any $r\in[k-1]$.

In order to show that $\pi(k)=1$ is not \esp-reducible, consider the gap vector $\g=(0,0,\dots,0,3,0)$.  Then $|Z(\{1234,4231\};\pi;\g)|<|Z(\{1234,4231\};\d_k(\pi);\d_k(\g))|$, finishing the proof.
\end{proof}

The substitution decomposition approach appears, at least on the surface, completely independent from the other three methods.  The class $\Av(321,2341,3412,4123)$ has a finitely labeled generating tree by Theorem~\ref{finlabel} (and thus it also has a finite enumeration scheme and a regular insertion encoding), but it contains infinitely many simple permutations\footnote{These simple permutations can be defined as the standardizations of even-length initial segments of the sequence $4,1,6,3,8,5,\dots,2k+2,2k-1,\dots$.}.

Moreover, the separable permutations defined in Section~\ref{wp-simple}, which contain only three simple permutations, possess neither a finite enumeration scheme nor a regular insertion encoding.  To see that they do not possess a regular insertion encoding, we need only note that they have an algebraic generating function (so, for this purpose, $\Av(132)$ would work just as well).

To establish that this class does not have a finite enumeration scheme, thereby completing our list of negative examples, we show that every permutation which consists of an increasing sequence followed by a decreasing sequence is \esp-irreducible for $\Av(2413,3142)$.  (This set forms a permutation class itself, with basis $\{213,312\}$.)  Therefore, not only do the separable permutations not have a finite enumeration scheme, but for each $n$ there are $2^{n-1}$ \esp-irreducible permutations.  Moreover, one can observe (either from the definition or the basis) that this class is invariant under the eight permutation class symmetries, so none of these offer any simplification.

\begin{proposition}\label{wp:serparable:bad}
Every $\pi\in \Av_k(213,312)$ is \esp-irreducible for $\Av(2413,3142)$.
\end{proposition}
\begin{proof}
Take $\pi\in\Av_k(213,312)$.  First we show that $\pi(r)$ is not \esp-reducible for any $2\le r\le k$.  If $\pi(r-1)>\pi(r)$, consider the gap vector $\g$ with all components $0$ except for $g_{r}=g_{r+1}=1$.  Then $Z(\{2413,3142\};\pi;\g)$ contains at most one permutation, because to avoid $2413$, the smaller of these new entries must be to the left of the larger one.  However, $Z(\{2413,3142\};d_r(\pi);d_r(\g))$ contains $2$ permutations (the most possible).  The case where $\pi(r-1)<\pi(r)$ can be handled similarly with the gap vector that is $0$ except for $g_{r-1}=g_r=1$.  Finally, $\pi(1)$ is not \esp-reducible because $Z(\{2413,3142\};\pi;(2,1,0,\dots,0))$ contains $5$ permutations -- every ordering of the new entries is allowed except $132$, because that would give rise to a $2413$-pattern -- while $Z(\{2413,3142\};d_1(\pi);(3,0,\dots,0))$ contains $6$ permutations.
\end{proof}

\section{An assortment of enumeration schemes successes}\label{wp-ex}

\begin{figure}[ht!]
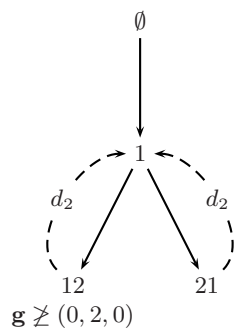

\begin{footnotesize}
\begin{center}
\pstree[nodesep=3pt,treefit=loose,treesep=10pt,levelsep=50pt]{\TR[name=0]{$\emptyset$}}{
\pstree{\TR[name=1,edge={\ncline{->}}]{$1$}}{
	\TR[name=12,edge={\ncline{->}}]{\basisnodeone{12}{(0,2,0)}}
	\TR[name=21,edge={\ncline{->}}]{$21$}
\nccurve[linestyle=dashed,angleA=135,angleB=180]{->}{12}{1}
\ncput*{$d_2$}
\nccurve[linestyle=dashed,angleA=45,angleB=0]{->}{21}{1}
\ncput*{$d_2$}}}
\end{center}
\end{footnotesize}
\caption{The enumeration scheme for $\Av(1342,1432)$}\label{schroeder-fig}
\end{figure}

Here we present several examples of finite enumeration schemes.  In these presentations, we adopt the following pictorial representation.  If $\pi$ is \esp-irreducible, so $|Z(B;\pi;\g)|$ is computed by summing over the $B$-avoiding children of $\pi$, then we draw a solid arrow from $\pi$ to each of its children.  If the entry $\pi(r)$ if \esp-reducible in $\pi$ then we draw a dashed arrow from $\pi$ to $d_r(\pi)$, label this arrow with $d_r$, and indicate the basis of $G_r(\pi)$ beneath $\pi$.  For example, the enumeration scheme for $\Av(1342,1432)$ is shown in Figure~\ref{schroeder-fig}.  We also define the {\it depth\/} of an enumeration scheme to be the least integer $k$ so that every permutation of length at least $k$ is \esp-reducible.

\begin{figure}[ht!]
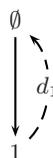

\begin{footnotesize}
\begin{center}
\pstree[nodesep=3pt,treefit=loose,treesep=10pt,levelsep=50pt]{\TR[name=0]{$\emptyset$}}{
	\TR[name=1,edge={\ncline{->}}]{$1$}
\nccurve[linestyle=dashed,angleA=45,angleB=315]{->}{1}{0}
\ncput*{$d_1$}}
\end{center}
\end{footnotesize}
\caption{The enumeration scheme for the set of all permutations, $\Av(\emptyset)$}\label{all-perms-fig}
\end{figure}

\begin{figure}[ht!]
\begin{footnotesize}
\begin{center}
\pstree[nodesep=3pt,treefit=loose,treesep=10pt,levelsep=50pt]{\TR[name=0]{$\emptyset$}}{
	\TR[name=1,edge={\ncline{->}}]{\basisnodeone{1}{(a,b)}}
\nccurve[linestyle=dashed,angleA=45,angleB=315]{->}{1}{0}
\ncput*{$d_1$}}
\end{center}
\end{footnotesize}
\caption{The enumeration scheme for $\Av(M_{a,b})$}\label{mansour-1-fig}
\end{figure}

As can be seen in Figure~\ref{schroeder-fig}, we include the empty permutation, $\emptyset$, in our diagrams.  Although this rarely has no more effect than making our diagrams consume more vertically space, there are a few exceptions.  One is the set of all permutations, $\Av(\emptyset)$, which has the enumeration shown in Figure~\ref{all-perms-fig}.  Another exception is $\Av(M_{a,b})$ where $M_{a,b}=\{\beta\in S_{a+b+1} : \beta(a+1)=1\}$.  This class was first counted by Mansour~\cite{m:cont:and:avoid}.  Its enumeration scheme is shown in Figure~\ref{mansour-1-fig}.

The scheme for $\Av(1234)$, which generates sequence \OEISlink{A005802} in the \href{http://www.research.att.com/\~njas/sequences/}{OEIS}~\cite{OEIS}, is shown in Figure~\ref{1234-fig}.

The $321$, hexagon-avoiding permutations, which can be defined as
$$
\Av(321,46718235,46781235,56718234, 56781234),
$$
were first introduced by Billey and Warrington~\cite{bw:321hex}, who showed how to compute the Kazhdan-Lusztig polynomials for them.  Stankova and West~\cite{sw:321hex} proved that the number, $s_n$, of these permutations of length $n$ satisfies
$$
s_n=6s_{n-1}-11s_{n-2}+9s_{n-3}-4s_{n-4}-4s_{n-5}+s_{n-6}
$$
for all $n\ge 7$, which gives sequence \OEISlink{A058094} in the \OEISref.  (So this complicated scheme produces only a sequence with a rational generating function.)  Later, Mansour and Stankova~\cite{ms:321hex} counted $321$, $2k$-gon-avoiding permutations for all $k$.

\begin{figure}[t]
\begin{footnotesize}
\pstree[nodesep=3pt,treefit=loose,treesep=10pt,levelsep=50pt]{\TR[name=0]{$\emptyset$}}{
\pstree{\TR[name=1,edge={\ncline{->}}]{$1$}}{
	\pstree{\TR[name=12,edge={\ncline{->}}]{$12$}}{
		\TR[name=123,edge={\ncline{->}}]{\basisnodeone{123}{(0,0,0,1)}}
		\TR[name=132,edge={\ncline{->}}]{$132$}
		\pstree{\TR[name=312,edge={\ncline{->}}]{$312$}}
			{%
			\TR[name=3124,edge={\ncline{->}}]{\basisnodeone{3124}{(0,0,0,0,1)}}
			\TR[name=4312,edge={\ncline{->}}]{$4312$}
			\TR[name=3142,edge={\ncline{->}}]{$3142$}
			\TR[name=3412,edge={\ncline{->}}]{$3412$}
			}
	}
	{%
	\pstree{\TR[name=21,edge={\ncline{->}}]{$21$}}
		{
		\TR[name=231,edge={\ncline[linestyle=none]}]{$231$}
		}
	}
\nccurve[linestyle=dashed,angleA=90,angleB=0]{->}{21}{1}
\ncput*{$d_2$}
\nccurve[linestyle=dashed,angleA=90,angleB=180]{->}{123}{12}
\ncput*{$d_3$}
\nccurve[linestyle=dashed,angleA=0,angleB=270]{->}{132}{12}
\ncput*{$d_3$}
\nccurve[linestyle=dashed,angleA=90,angleB=270]{->}{3412}{231}
\ncput*{$d_3$}
\nccurve[linestyle=dashed,angleA=90,angleB=180]{->}{3124}{312}
\ncput*{$d_4$}
\nccurve[linestyle=dashed,angleA=90,angleB=180]{->}{3142}{231}
\ncput*{$d_2$}
\nccurve[linestyle=dashed,angleA=150,angleB=255]{->}{4312}{312}
\ncput*{$d_2$}
\nccurve[linestyle=dashed,angleA=150,angleB=0]{->}{231}{12}
\ncput*{$d_3$}
}}
\end{footnotesize}
\caption{The enumeration scheme for $\Av(1234)$}\label{1234-fig}
\end{figure}

\begin{figure}[t]
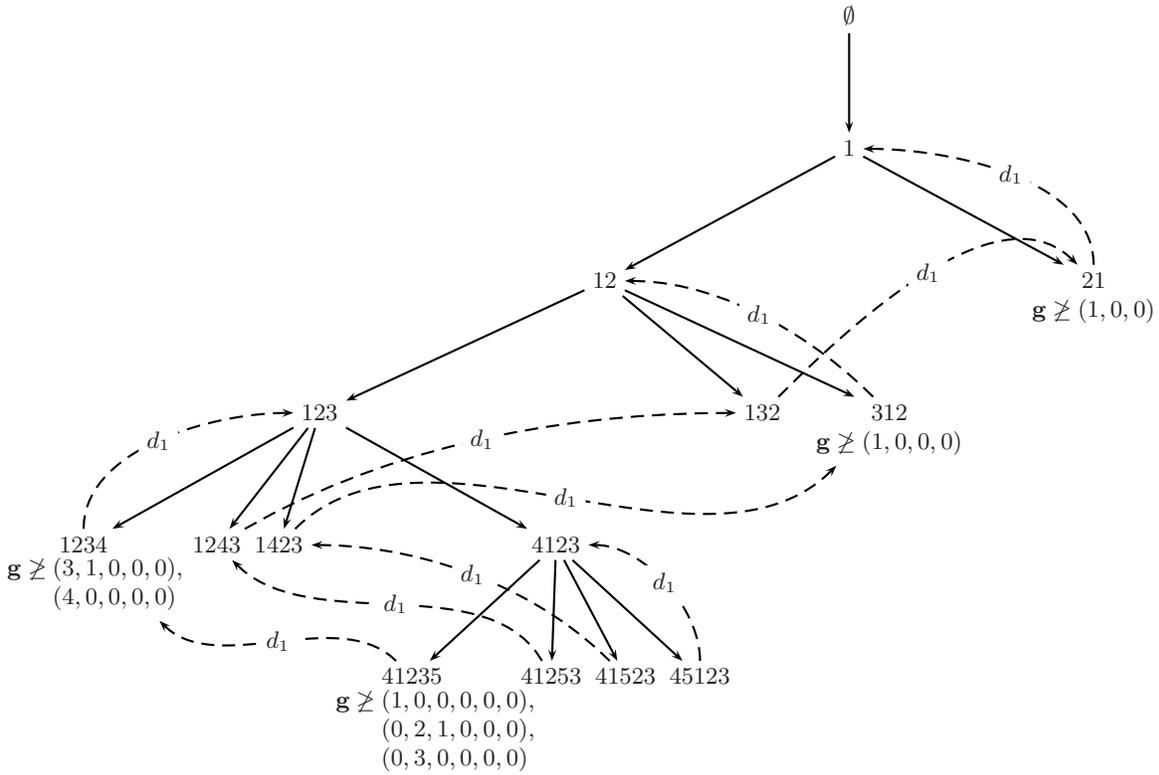

\begin{footnotesize}
\pstree[nodesep=3pt,treefit=loose,treesep=5pt,levelsep=50pt]{\TR[name=0]{$\emptyset$}}{
\pstree{\TR[name=1,edge={\ncline{->}}]{$1$}}{
	\pstree{\TR[name=12,edge={\ncline{->}}]{$12$}}
	{%
		\pstree{\TR[name=123,edge={\ncline{->}}]{$123$}}
		{%
			\TR[name=1234,edge={\ncline{->}}]{\basisnodetwo{1234}{(3,1,0,0,0)}{(4,0,0,0,0)}}
			\TR[name=1243,edge={\ncline{->}}]{$1243$}
			\TR[name=1423,edge={\ncline{->}}]{$1423$}
			\pstree{\TR[name=4123,edge={\ncline{->}}]{$4123$}}
			{%
				\TR[name=41235,edge={\ncline{->}}]{%
				  	\basisnodethree{41235}{(1,0,0,0,0,0)}{(0,2,1,0,0,0)}{(0,3,0,0,0,0)}}
				\TR[name=41253,edge={\ncline{->}}]{$41253$}
				\TR[name=41523,edge={\ncline{->}}]{$41523$}
				\TR[name=45123,edge={\ncline{->}}]{$45123$}
			}
		}
		\TR[name=132,edge={\ncline{->}}]{$132$}
		\TR[name=312,edge={\ncline{->}}]{\basisnodeone{312}{(1,0,0,0)}}
	}
	\TR[name=21,edge={\ncline{->}}]{\basisnodeone{21}{(1,0,0)}}
	\nccurve[linestyle=dashed,angleA=90,angleB=0]{->}{21}{1}
	\ncput*{$d_1$}
	\nccurve[linestyle=dashed,angleA=45,angleB=135]{->}{132}{21}
	\ncput*{$d_1$}
	\nccurve[linestyle=dashed,angleA=135,angleB=0]{->}{312}{12}
	\ncput*{$d_1$}
	\nccurve[linestyle=dashed,angleA=90,angleB=180]{->}{1234}{123}
	\ncput*{$d_1$}
	\nccurve[linestyle=dashed,angleA=30,angleB=180]{->}{1243}{132}
	\ncput*{$d_1$}
	\nccurve[linestyle=dashed,angleA=45,angleB=225]{->}{1423}{312}
	\ncput*{$d_1$}
	\nccurve[linestyle=dashed,angleA=135,angleB=315]{->}{41235}{1234}
	\ncput*{$d_1$}
	\nccurve[linestyle=dashed,angleA=120,angleB=315]{->}{41253}{1243}
	\ncput*{$d_1$}
	\nccurve[linestyle=dashed,angleA=135,angleB=0]{->}{41523}{1423}
	\ncput*{$d_1$}
	\nccurve[linestyle=dashed,angleA=90,angleB=0]{->}{45123}{4123}
	\ncput*{$d_1$}
}}
\end{footnotesize}
\caption{The enumeration scheme for the $321$, hexagon-avoiding permutations}\label{321-hex-fig}
\end{figure}

\begin{figure}[t]
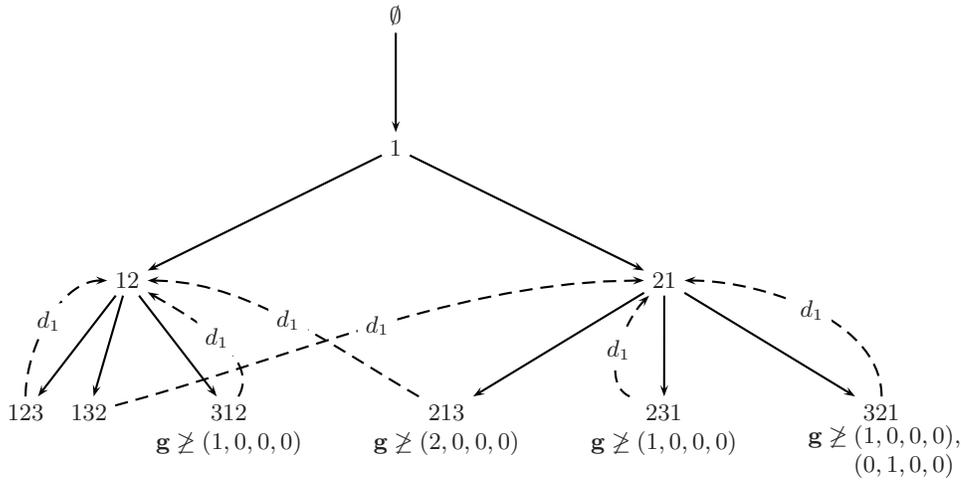

\begin{footnotesize}
\pstree[nodesep=3pt,treefit=loose,treesep=10pt,levelsep=50pt]{\TR[name=0]{$\emptyset$}}{
\pstree{\TR[name=1,edge={\ncline{->}}]{$1$}}{
\pstree{\TR[name=12,edge={\ncline{->}}]{$12$}}{
\TR[name=123,edge={\ncline{->}}]{$123$}
\TR[name=132,edge={\ncline{->}}]{$132$}
\TR[name=312,edge={\ncline{->}}]{\basisnodeone{312}{(1,0,0,0)}}
}
{\pstree{\TR[name=21,edge={\ncline{->}}]{$21$}}{
\TR[name=213,edge={\ncline{->}}]{\basisnodeone{213}{(2,0,0,0)}}
\TR[name=231,edge={\ncline{->}}]{\basisnodeone{231}{(1,0,0,0)}}
\TR[name=321,edge={\ncline{->}}]{\basisnodetwo{321}{(1,0,0,0)}{(0,1,0,0)}}
}}}
\nccurve[linestyle=dashed,angleA=90,angleB=0]{->}{321}{21}
\ncput*{$d_1$}
\nccurve[linestyle=dashed,angleA=155,angleB=225]{->}{231}{21}
\ncput*{$d_1$}
\nccurve[linestyle=dashed,angleA=150,angleB=0]{->}{213}{12}
\ncput*{$d_1$}
\nccurve[linestyle=dashed,angleA=60,angleB=330]{->}{312}{12}
\ncput*{$d_1$}
\nccurve[linestyle=dashed,angleA=15,angleB=180]{->}{132}{21}
\ncput*{$d_1$}
\nccurve[linestyle=dashed,angleA=90,angleB=180]{->}{123}{12}
\ncput*{$d_1$}
}
\end{footnotesize}
\caption{The enumeration scheme for the freely braided permutations, $\Av(3421,4231,4312,4321)$}
\label{freely-braided-fig}
\end{figure}

The {\it freely braided permutations\/} are the class
$$
\Av(3421,4231,4312,4321).
$$
This class was introduced by Green and Losonczy~\cite{gl:fb} and also arises in the work of Tenner~\cite{ten:rdp}.  Mansour~\cite{mansour:fb} found that the freely braided permutations have the generating function
$$
\frac{1-3x-2x^2+(1+x)\sqrt{1-4x}}{1-4x-x^2+(1-x^2)\sqrt{1-4x}}.
$$
\OEIS{A108600}  They also have an enumeration scheme of depth $3$, shown in Figure~\ref{freely-braided-fig}.

\begin{figure}[t]
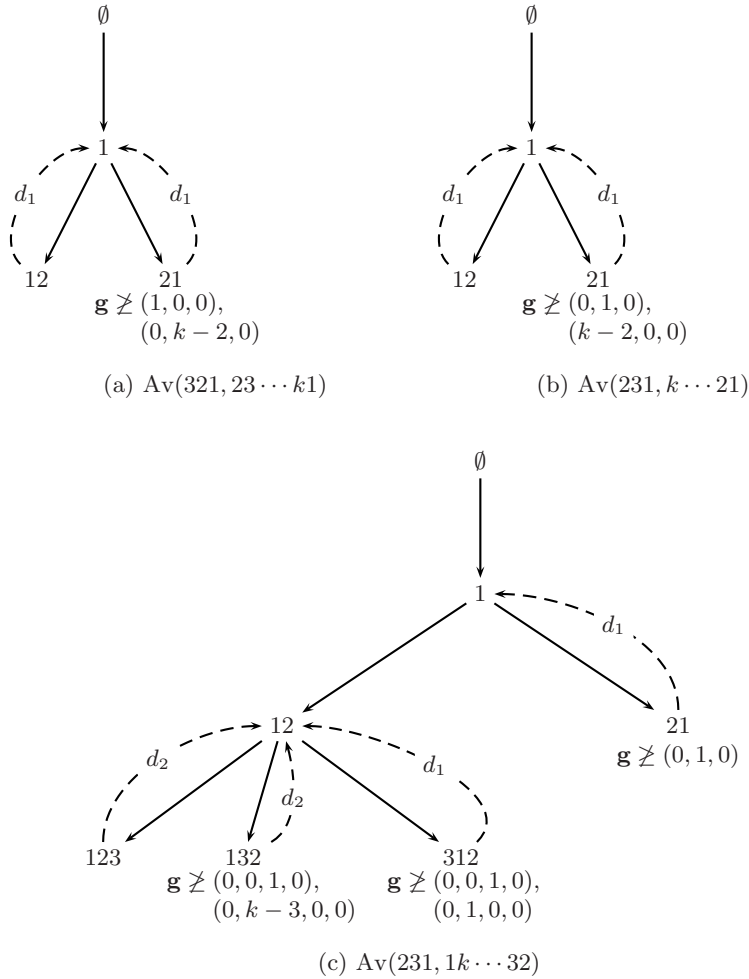

\begin{tabular}{cc}
\subfigure[$\Av(321,23\cdots k1)$]{
	\parbox{2in}{
		\begin{footnotesize}
		\pstree[nodesep=3pt,treefit=loose,treesep=10pt,levelsep=50pt]{
		\TR[name=0]{$\emptyset$}}{
		\pstree{\TR[name=1,edge={\ncline{->}}]{$1$}}{
			\TR[name=12,edge={\ncline{->}}]{$12$}
			\TR[name=21,edge={\ncline{->}}]{\basisnodetwo{21}{(1,0,0)}{(0,k-2,0)}}
		\nccurve[linestyle=dashed,angleA=135,angleB=180]{->}{12}{1}
		\ncput*{$d_1$}
		\nccurve[linestyle=dashed,angleA=45,angleB=0]{->}{21}{1}
		\ncput*{$d_1$}}}
		\end{footnotesize}
	}
}
&
\subfigure[$\Av(231,k\cdots 21)$]{
	\parbox{2in}{
		\begin{footnotesize}
		\pstree[nodesep=3pt,treefit=loose,treesep=10pt,levelsep=50pt]{
		\TR[name=0]{$\emptyset$}}{
		\pstree{\TR[name=1,edge={\ncline{->}}]{$1$}}{
			\TR[name=12,edge={\ncline{->}}]{12}
			\TR[name=21,edge={\ncline{->}}]{\basisnodetwo{21}{(0,1,0)}{(k-2,0,0)}}
		\nccurve[linestyle=dashed,angleA=135,angleB=180]{->}{12}{1}
		\ncput*{$d_1$}
		\nccurve[linestyle=dashed,angleA=45,angleB=0]{->}{21}{1}
		\ncput*{$d_1$}}}
		\end{footnotesize}
	}
}
\\
\multicolumn{2}{c}{
\subfigure[$\Av(231,1k\cdots 32)$]{
	\begin{footnotesize}
	\pstree[nodesep=3pt,treefit=loose,treesep=10pt,levelsep=50pt]{
	\TR[name=0]{$\emptyset$}}{
	\pstree{\TR[name=1,edge={\ncline{->}}]{$1$}}{
		\pstree{\TR[name=12,edge={\ncline{->}}]{$12$}}{
			\TR[name=123,edge={\ncline{->}}]{123}
			\TR[name=132,edge={\ncline{->}}]{\basisnodetwo{132}{(0,0,1,0)}{(0,k-3,0,0)}}
			\TR[name=312,edge={\ncline{->}}]{\basisnodetwo{312}{(0,0,1,0)}{(0,1,0,0)}}
		}
		\TR[name=21,edge={\ncline{->}}]{\basisnodeone{21}{(0,1,0)}}
	\nccurve[linestyle=dashed,angleA=90,angleB=0]{->}{21}{1}
	\ncput*{$d_1$}
	\nccurve[linestyle=dashed,angleA=90,angleB=180]{->}{123}{12}
	\ncput*{$d_2$}
	\nccurve[linestyle=dashed,angleA=30,angleB=285]{->}{132}{12}
	\ncput*{$d_2$}
	\nccurve[linestyle=dashed,angleA=45,angleB=0]{->}{312}{12}
	\ncput*{$d_1$}}}
	\end{footnotesize}
}}
\end{tabular}
\caption{Enumeration schemes for three Wilf-equivalent classes}\label{chow-west-fig}
\end{figure}

Chow and West~\cite{cw:cheby} showed using generating trees that the classes in Figure~\ref{chow-west-fig} are Wilf-equivalent, that their generating functions are rational\footnote{This fact can be verified quite quickly: by inverting the permutations in $\Av(321,23\cdots k1)$ one obtains the class $\Av(321,k12\cdots (k-1))$, which satisfies the hypotheses of Theorem~\ref{finlabel}, and thus also of Theorem~\ref{insertion}.  Note, however, that while symmetries of these classes always have finitely labeled generating trees, and thus also regular insertion encodings, the complexity of their generating trees and insertion encodings increases with $k$, while the complexity of their enumeration schemes stay fixed.}, and that these generating functions can be expressed as quotients of Chebyshev polynomials.  Krattenthaler~\cite{kratt:cheby} gave another proof via a bijection to Dyck paths (in fact, he found the bivariate generating functions for $132$-avoiding permutations by length and number of copies of $12\cdots k$ and for $123$-avoiding permutations by length and number of copies of $(k-1)\cdots 21k$).  Around the same time as that work, several other authors studied these and similar classes: Jani and Rieper~\cite{jr:catalan}, Mansour and Vainshtein~\cite{mv:cheby1, mv:cheby2}, and Robertson, Wilf, and Zeilberger~\cite{rwz:fractions}.  Deutsch, Hildebrand, and Wilf~\cite{dhw:lis} used these results in their study of the longest increasing subsequence problem for $132$-avoiding permutations.

For any fixed $k$, {\sc WilfPlus} can automatically and rigorously derive enumeration schemes for these classes.  This makes it easy to make conjectures for the general form of these enumeration schemes, although they must be verified by hand.  We carry this out for one of these classes below.

\begin{proposition}\label{prop-chow-west}
The enumeration scheme for $\Av(231,k\cdots 21)$ is as shown in Figure~\ref{chow-west-fig} (b).
\end{proposition}
\begin{proof}
Let $B=\{231,k\cdots 21\}$.  First we verify the computations of $G_r(\pi)$ given in the diagram.  For $12$, the diagram shows that $\G_1(12)=\mathbb{N}^3$.  In any permutation in $Z(B;12;(g_1,g_2,g_3))$, the entries before the $2$ must be in decreasing order from left to right as any ascent before the $2$ would give rise to a $231$-pattern.  Then, since these permutations must avoid $k\cdots 21$, there can be at most $k-2$ new entries before the $2$, so $(g_1,g_2,g_3)\notin \G(12)$ whenever $g_1+g_2\ge k-1$.  The $B$-avoiding permutations
$$
g_1+g_2+2,\dots,g_2+4,g_2+3,1,g_2+2,\dots,4,3,2,g_1+g_2+3,g_1+g_2+4,\dots,g_1+g_2+g_3+2
$$
then show that $\G(12)$ is $\Av(\{(g_1,g_2,g_3) : g_1+g_2\ge k-1\})$.  It can be shown similarly that $\G(1)=\Av((k-1,0))$.  This shows that $\G_1(12)=\mathbb{N}^3$ because $Z(B;1;\d_1(\g))$ is empty whenever $Z(B;12;\g)$ is empty.

The computation for $21$ is similar.  One first checks that
$$
\G(21)=\Av((0,1,0),(k-2,0,0)).
$$
Since $\G(1)=\Av((k-1,0))$, the set of gap vectors $\g$ for which either $Z(B;21;\g)\neq\emptyset$ or $Z(B;1;\d_1(\g))=\emptyset$ is $\Av((0,1,0))\setminus\{(k-2,0,0)\}$.  This implies that $\G_1(21)=\Av((0,1,0),(k-2,0,0))$.
Now we must verify the \esp-reducible entries.  The $12$ case is clear because the only way $1$ can participate in a $231$ or $k\cdots 21$-pattern is as the minimal entry, and in either case the $2$ could play the same role.  For $21$, first note that the $2$ also cannot participate in any $231$ or $k\cdots 21$-pattern if $\g\in\G_1(21)$: it cannot be the ``$2$'' in a $231$ because there can be nothing between the $2$ and the $1$, and it cannot be the ``$2$'' in a $k\cdots 21$-pattern because our restrictions on $\g$ prevent there from being enough entries to the left of the $2$ to accommodate such a pattern.  This shows that the $2$ can only possibly be the minimal entry in a forbidden pattern, but then the $1$ could play the same role.
\end{proof}

In addition to the examples already mentioned, there are several other interesting classes avoiding a pair of permutations of length four\footnote{The characterization of the Wilf-equivalences between these classes has recently been completed by Le~\cite{le:length4}.}:
\begin{enumerate}
\item Kremer and Shiu~\cite{ks:len4} proved that there are four classes of this form counted by $(4^{n-1}+2)/3$.  These all have finite enumeration schemes: $\Av(1234,2143)$ and $\Av(1432,2341)$ have schemes of depth 3, $\Av(2341,4321)$ has a scheme of depth 4, and $\Av(2143,4123)$ has a scheme of depth 6.
\item Kremer~\cite{k:sn} (see also Stanley~\cite[Exercise 6.39.l]{stanley:ec2}) completed the characterization of the classes defined by avoiding two permutations of length four that are counted by the large Schr\"oder numbers.  Up to symmetry, there are ten such classes, seven of which have finite enumeration schemes:
\medskip
\begin{footnotesize}
\begin{center}
\begin{tabular}{l|l}
Class&Finite enumeration scheme?\\
\hline\hline
$\Av(1342,2341)$&Yes, of depth 3\\
$\Av(1342,1432)$&
\begin{minipage}[t]{2.2in}
Yes, of depth 2, shown in Figure~\ref{schroeder-fig}
\end{minipage}\\
$\Av(2341,2413)$&No\\
$\Av(2413,3142)$&
\begin{minipage}[t]{2.2in}
No, these are the separable permutations considered in Proposition~\ref{wp:serparable:bad}
\end{minipage}
\\
$\Av(2431,3241)$&No\\
$\Av(3241,3421)$&Yes, of depth 4\\
$\Av(3241,4231)$&
\begin{minipage}[t]{2.2in}
Yes, of depth 2 (Knuth~\cite{knuth1} proved that these are precisely the permutations that can be sorted by an input-restricted deque)
\end{minipage}
\\
$\Av(3412,3421)$&Yes, of depth 3\\
$\Av(3421,4321)$&Yes, of depth 2\\
$\Av(3421,4231)$&Yes, of depth 4
\end{tabular}
\end{center}
\end{footnotesize}
\item A permutation is said to be {\it skew-merged\/} if it is the union of an increasing subsequence and a decreasing subsequence.  Stankova~\cite{stankova:fs} was the first to prove that the skew-merged permutations are given by $\Av(2143,3412)$.  Later K\'ezdy, Snevily, and Wang~\cite{ksw:incdec} gave another proof of this result using F\"oldes and Hammer's characterization of split graphs~\cite{fh:sg}.  Atkinson~\cite{a:skewmerged} showed that the generating function for this class is
$$
\frac{1-3x}{(1-2x)\sqrt{1-4x}},
$$
(sequence \OEISlink{A029759} in the \href{http://www.research.att.com/\~njas/sequences/}{OEIS}~\cite{OEIS}).  This class has an enumeration scheme of depth four.
\end{enumerate}

\section{Conclusion}\label{wp-conclusion}

We have developed an extension of Zeilberger's enumeration schemes and provided the mechanisms for their automatic generation, a task implemented in the accompanying Maple package {\sc  WilfPlus}, which is available at \url{http://math.rutgers.edu/~vatter/}.  As the examples in Section~\ref{wp-nonex} demonstrate, there remain considerable differences in the applicability of enumeration schemes, substitution decompositions, and the insertion encoding.  This obviously suggests the following question.

\begin{question}\label{systematic-gen}
Is there a systematic method of permutation class enumeration which is applicable to all classes with finite enumeration schemes, all classes with regular insertion encodings, and all classes with only finitely many simple permutations?
\end{question}

Additionally, one would like such a method to be invariant under the eight permutation class symmetries, but Question~\ref{systematic-gen} is probably demanding enough in the form above.

We conclude by collecting three of the questions raised earlier.

\begin{question}\label{prec}
Is every sequence produced by a finite enumeration scheme holonomic?
\end{question}

As demonstrated by the examples in Section~\ref{wp-ex}, another interesting question is the equivalence problem:

\begin{question}
Is it decidable whether two finite enumeration schemes produce the same sequence?
\end{question}

Brlek, Duchi, Pergola, and Rinaldi~\cite{bdpr:eps} consider the equivalence problem for generating trees with infinitely many labels.

Although only tangentially related to the main thrust of this article, the following question nevertheless seems intriguing enough to warrant its inclusion.

\begin{question}
Is it decidable whether a class contains only finitely many simple permutations?
\end{question}

\bigskip
\bibliographystyle{acm}
\bibliography{wilfplus}

\end{document}